\documentclass[11pt]{amsart}

\usepackage{amsmath,amssymb,amsthm,graphicx}

\usepackage{geometry}
\usepackage{color}
\theoremstyle{definition}
\newtheorem{defn}{Definition}[section]

\theoremstyle{plain}
\newtheorem{theorem}[defn]{Theorem}

\newcommand{\qb}[2]{\left[\begin{array}{c} #1 \\ #2  \end{array}\right]}
\newcommand{\s}{\sigma}
\newcommand{\Ga}{\Gamma}

\begin{document}
\title{Colored Jones polynomials without tails}

\author[C. Lee]{Christine Ruey Shan Lee}
\author[R. van der Veen]{Roland van der Veen}

\address[]{Department of Mathematics, University of Texas at Austin, Austin TX 78712}
\email[]{clee@math.utexas.edu}

\address[]{Mathematisch Instituut, Leiden University, Leiden, Netherlands}
\email[]{r.i.van.der.veen@math.leidenuniv.nl}

\thanks{Lee was partially supported by NSF Grant MSPRF-DMS 1502860. Van der Veen was supported by the Netherlands foundation for scientific research (NWO)}

\begin{abstract}We exhibit an infinite family of knots with the property that the first coefficient of the $n$-colored Jones polynomial grows linearly with $n$.
This shows that the concept of stability and tail seen in the colored Jones polynomials of alternating knots does not generalize naively.
\end{abstract}

\maketitle

\section{Introduction}
%Will ask Abhijit about how to cite their results.
The colored Jones polynomials are a well known knot invariant with many connections to other parts of mathematics \cite{Jo,Tu}.
It is the Reshetikhin-Turaev knot invariant corresponding to the $n$-dimensional representation of quantum $sl_2$.
Given a knot $K$ and a natural number $n\geq 2$ the $n$-colored Jones polynomial is a Laurent polynomial $J_{K,n}(v)$. Dasbach-Lin \cite{DL06} observed that for any fixed alternating knot the first coefficients seem to stabilize as $n$ grows. The sequence of stabilized coefficients was called the tail of the knot. Tails are defined and studied in greater detail in \cite{Arm13,CK13,GL15, AD17, Haj17, L18}.

For general knots similar behaviour was expected\footnote{See for example Remark 1.7 of \cite{GaVu17}.} provided one allowed a certain periodicity. For example the leading coefficient of most torus knots is $(-1)^{n}$, and they are shown to have one or more tails \cite{AD18, GL15}; for an index $b$ closed braid consisting of a positive full twist on $b$ strands and $\ell\leq b$ negative crossings, the work of \cite{CK13} determines that up to $b-\ell+1$ first coefficients of its colored Jones polynomial are stable, though it might still admit multiple tails. 

\begin{figure}
\begin{center}
\includegraphics[width=5cm]{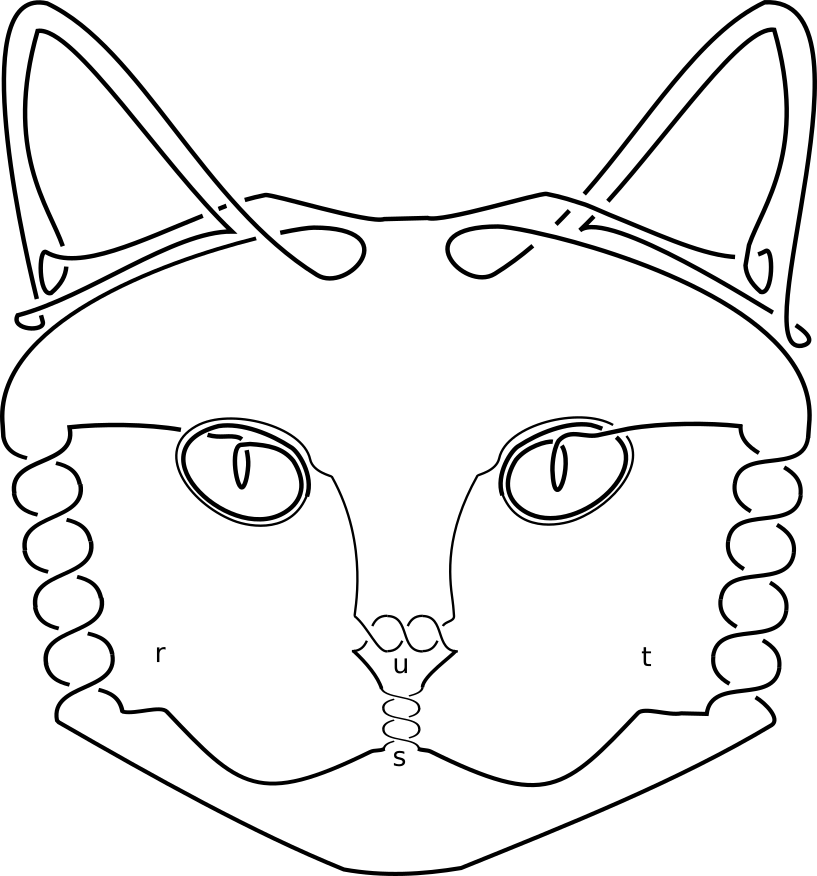}
\end{center}
\caption{The family of knots $C(r,s,t,u)$ where $r,s,t,u$ denote the number of crossings in each twist region. We have drawn the Manx knot $C(5,-3,5,-3)$.
To make the diagram zero framed we drew $16=r-s+t-u$ additional negative curls.}
\label{fig.ManxKnot}
\end{figure}

In this note we introduce knots where the first coefficient of the colored Jones polynomial does not stabilize at all. Instead it grows without bound, proportional to $n$. Such knots can never have tails, so we propose to name such knots Manx knots for a breed of cats without tails.
%Should we say that it is the simplest Manx pretzel? 
Possibly the simplest Manx knot is the pretzel knot $P(5,-3,5)$. The infinite family of Manx knots $C(5,-3,5,u)$ introduced in this paper are the Montesinos knots shown in Figure \ref{fig.ManxKnot}.

We will first give a quick, self-contained definition of the colored Jones polynomial that allows us to compute a formula for the polynomial of our examples of Manx knots. Next we will analyse the degree of the terms in our formula and observe the lack of tail in the main theorem, Theorem \ref{t.main}.

\section{Colored Jones polynomial}

The colored Jones polynomial is usually considered as a knot invariant but from our point of view it is more natural to consider its extension to knotted framed graphs.
A good thing about this definition of the colored Jones polynomial is that it is less dependent on knot diagrams and a little more efficient for computations. The technical reason behind its success is that we diagonalize the R-matrix.

\begin{defn}Knotted trivalent graphs (KTG)
\begin{enumerate}	
\item A framed graph is a 1-dimensional simplicial complex $\Gamma$ together with an embedding $\Gamma\to \Sigma$ of $\Gamma$ into a surface with boundary $\Sigma$ as a spine.
\item A coloring of $\Gamma$ is a map $\sigma:E(\Gamma)\to \mathbb{N}$.
\item A Knotted Trivalent Graph (KTG) is a trivalent framed graph embedded (as a surface) into $\mathbb{R}^3$, considered up to isotopy.
\end{enumerate}
\end{defn}

Framed knots and links are special cases of KTGs. We prefer the more general set of KTGs here because there are natural 3-dimensional operations
that can be used to generate the whole theory from the tetrahedron. 

\begin{figure}[htp]
\begin{center}
\def \svgwidth{.8\columnwidth}
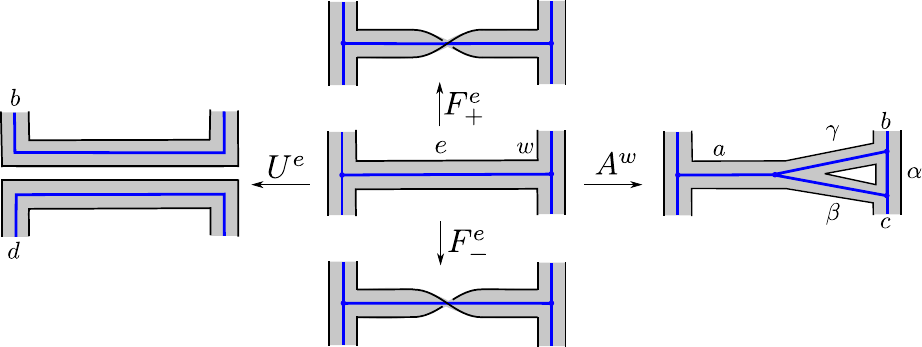
\end{center}
\caption{Operations on Knotted Trivalent Graphs: framing change $F^e_{\pm}$, unzip $U^e$ and triangle move $A^w$ applied to an edge $e$ and vertex $w$ shown in the middle.}
\label{fig.KTGMoves}
\end{figure}

More precisely it was shown by D. Thurston \cite{ThuKTG} and \cite{Vdv} that the following 
three moves suffice to generate any KTG from the planar theta graph. The first operation $F^e_{\pm}$ cuts the surface $\Sigma$ transversal to an edge, rotates one side by $\pi$
and reglues. This is called framing change or connected sum with a M\"{o}bius band. The second operation is called unzip, $U^e$. It doubles a chosen edge along its framing, deletes its end-vertices and joins the result as shown. The final move is called $A^w$ and expands a vertex $w$ into a triangle. All three moves are as shown in Figure \ref{fig.KTGMoves}. 

Define $[k] =\frac{v^{2k}-v^{-2k}}{v^2-v^{-2}}$ and $[k]! = [1][2]\dots[k]$ for $k\in \mathbb{N}$ and $[k]! = 0$ if $k\notin \mathbb{N}$. Let
\[\qb{a_1+a_2+\dots+a_r}{a_1,a_2,\dots, a_r} = \frac{[a_1+a_2+\dots+a_r]!}{[a_1]!\dots[a_r]!}\]

The building blocks of the Jones polynomial are the following Laurent polynomials. For each we explicitly state the term of maximal degree for later use.
\[
f(a) = i^{-a}v^{\frac{-a(a+2)}{2}} \qquad O^k = (-1)^k[k+1] = (-1)^kv^{2k}+lot
\]
\[
\theta(a,b,c) = O^{\frac{a+b+c}{2}}\qb{\frac{a+b+c}{2}}{\frac{-a+b+c}{2},\frac{a-b+c}{2},\frac{a+b-c}{2}} = (-1)^{\frac{a+b+c}{2}}
v^{a+b+c+ab+bc+ca-\frac{a^2+b^2+c^2}{2}} +lot
\]
\[\Delta(a,b,c,\alpha,\beta,\gamma) = \sum_{z}\frac{O^z}{O^{\frac{a+b+c}{2}}} \qb{z}{\frac{a+b+c}{2}}\qb{\frac{-a+b+c}{2}}{z-\frac{a+\beta+\gamma}{2}}
\qb{\frac{a-b+c}{2}}{z-\frac{\alpha+b+\gamma}{2}}\qb{\frac{a+b-c}{2}}{z-\frac{\alpha+\beta+c}{2}} = \]
\[
(-1)^{m-\frac{a+b+c}{2}}v^{g(m+1,\frac{a+b+c}{2}+1)+g(\frac{-a+b+c}{2},m-\frac{a+\beta+\gamma}{2})+g(\frac{a-b+c}{2},m-\frac{\alpha+b+\gamma}{2})+g(\frac{a+b-c}{2},m-\frac{\alpha+\beta+c}{2})}+lot
\]
where $g(n,k) = 2k(n-k)$ and $2m = a+b+c+\alpha+\beta+\gamma-\max(a+\alpha,b+\beta,c+\gamma)$

Our formulas agree with the integer normalization of \cite{Cost}. The formula $\Delta$ is the quotient of the $6j$-symbol and a theta; it is in fact a Laurent polynomial. The summation range for $\Delta$ is finite as dictated by the binomials. 

The colored Jones polynomial of any KTG can now be defined by the following three equations. 
\begin{equation}
\langle F^e_{\pm}(\Ga),\s \rangle = f(\s(e))^{\pm1}\langle \Ga,\s \rangle
\end{equation}
Here $F_{\pm}^e(\Ga)$ means the result after applying framing change to edge $e$ of the KTG $\Ga$.
\begin{equation}
\langle U^e(\Ga),\s' \rangle = \langle \Gamma,\s \rangle\sum_{\s(e)}\frac{O^{\s(e)}}{\theta(\s(e),\s(b),\s(d))}
\end{equation}
The two colorings $\s$ and $\s'$ agree on all edges except that $U^e(\Ga)$ does not have edge $e$. The four edges in $\Ga$ adjacent to $e$ must have the same colors as $b$ and $d$ in $U^e(\Ga)$. 
Also, the summation is over all possible values of color of the missing edge $\s(e)$. Most values yield zero because the factorials vanish. The only values that may be non-zero are when $\s(e)$ is between $|\s'(b)-\s'(d)|$ and $\s'(b)+\s'(d)$ and has the same parity. Finally,
\begin{equation}
\langle A^w(\Ga),\s' \rangle = \langle \Ga,\s \rangle\Delta(\s'(a),\s'(b),\s'(c),\s'(\alpha),\s'(\beta),\s'(\gamma))
\end{equation}
In this equation the colorings $\s$ and $\s'$ again agree on all edges where it makes sense. The colors of the six new edges are denoted $a,\alpha,b,\beta,c,\gamma$ as shown in Figure \ref{fig.KTGMoves} and $\Delta$ 
is defined above. 
\begin{defn}
As noted above a $0$-framed knot $K$ is a special case of a $KTG$ with a single edge $e$. In this case we denote its colored Jones polynomial by 
$J_{K,n} = \langle K, \sigma \rangle$, where $\sigma(e)=n-1$.
\end{defn}

We study the family of Montesinos knots $C(r,s,t,u)$ shown in Figure \ref{fig.ManxKnot}. For technical reasons we restrict ourselves to $s\leq -3 \leq 3 r,t,u$ all odd. And $u<0$. A similar family was studied in \cite{LYL17}.

\begin{figure}[ht]
\tiny{
\begin{center}
\def \svgwidth{\columnwidth}
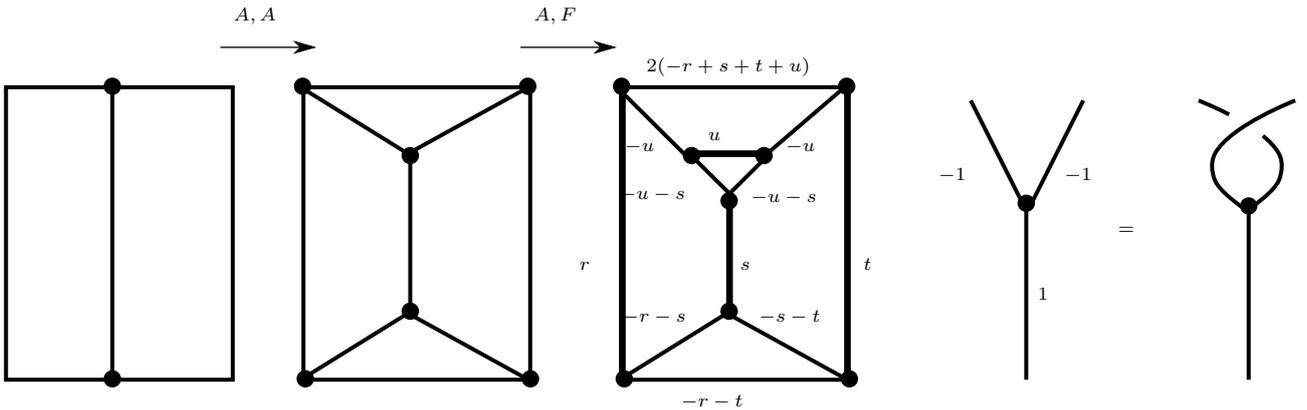
\end{center}}
\caption{\label{fig.KTGComputation} Starting from a theta graph (left), we first apply the $A$-move to both vertices, and again on one of the newly created vertices. Next change the framing on many edges. The numbers on each edge indicate the number of half twists present. Finally unzip the four edges (shown thickened) marked with $r, s, t, u$, respectively, to obtain the $0$-framed diagram for the knot $C(r,s,t,u)$.}
\end{figure}

To compute the colored Jones polynomial of $K = C(r,s,t,u)$ we take the following steps (see also \cite{LV16}). The first step is to generate our knot $K$ from the theta graph by KTG moves. One way to achieve this is shown in Figure \ref{fig.KTGComputation}. All edges should be thought of as bands with blackboard framing. The numbers indicate half twists (not colors). Note that the unzip applied to a twisted edge produces two twisted bands that form a crossing as shown in the figure on the right. This is natural considering that the black lines stand for actual strips.
Reading backwards and applying the above equations we may compute $J_{K,n+1}$ as follows. The unzips yield four summations, the framing change multiplies everything by the factors $f$, the $A$ moves produce the $\Delta$ terms and finally we are left with a single $\theta(2a,2b,2c)$. Note we have scaled the summation variables $a,b,c,d$ to make them integers instead of even integers.

\[
J_{C(r,s,t,u),n+1} = f(n)^{-2\mathrm{Wr}(r,s,t)-2r-2s-2t+2u}
\sum_{a,b,c,d \in nP\cap \mathbb{Z}}\Delta(2a,2b,2c,n,n,n)^2\times  \]
\[
\Delta(2b,n,n,2d,n,n)
\frac{O^{2a}O^{2b}O^{2c}O^{2d}
f(2a)^r f(2b)^s f(2c)^t f(2d)^{-u} 
\theta(2a,2b,2c)}{\theta(2a,n,n)\theta(2b,n,n)\theta(2c,n,n)\theta(2d,n,n)}
\]
Here the writhe is given by $\mathrm{Wr}(r,s,t,u)=r-s+t-u$ and \\$P = \{(a,b,c,d)\in [0,1]^4: |a-b|\leq c \leq a+b \}$.

\section{The lack of tail}

Using the formulas for the building blocks we may determine the degree $\delta(a,b,c,d)$ of each term in the state sum.
It is an explicit piecewise quadratic function on the polytope $nP$. From our assumptions it follows that:
%on $nP$ we have
\[\delta(a+1,b,c,d)<\delta(a,b,c,d) \text{ and } \delta(a,b,c+1,d)<\delta(a,b,c,d) < \delta(a,b,c,d-1).\] 
Therefore all terms of highest degree satisfy $b=a+c$ and $d = 0$ with $0\leq a,c \leq a+c\leq n$. Explicitly the degree is
$\delta(a,a+c,c,0) = $
\[b (r + 2 c (1 + r) + s-1)-2 (b^2 (1 + r + s)  + c (r + t-2) + 
    c^2 (r + t)) - 2 n (-2 r - u) - 2 n^2 (-r - u)\] %taken from CatKnots4.nb
Setting $r,s,t = 5,-3,5$ we find $\delta(a,b,c,d) = -4 (a - c)^2 + 2n(18 + 10 n - (2  + n) u)$
It follows that the maxima are on the line $a=c$ and there are $\lfloor n/2 \rfloor$ such maxima.
Looking at the building blocks, each of the corresponding maximal terms in the state sum has the same leading coefficient $(-1)^{8 a + 2 n (-11 + u)}=1$.
This completes the proof of our main theorem:

\begin{theorem} \label{t.main}
If $u\leq -3$ then $C(5,-3,5,u)$ is a Manx knot.\\
 More precisely, the leading coefficient of
$J_{C(5,-3,5,u),n+1}(v)$ is $\lfloor n/2 \rfloor$.
\end{theorem}

The same techniques will also show that the pretzel knot $P(5,-3,5)$ is Manx but this appears to be the only Manx pretzel knot \cite{GLV18}.
We expect that many more Manx knots exist but not in the family $C(r,s,t,u)$ we considered. All Montesinos Manx knots seem to derive from the unique Manx pretzel $P(5,-3,5)$.

Although our Manx knots do not stabilize in the ordinary sense they may still be $c$-stable. This more general notion of stability was introduced in Garoufalidis-Vuong \cite{GaVu17} and allows growing leading coefficients. Briefly, a sequence of Laurent power series $f_n(v)$ with $v$-degree $\delta(n)$ is $c$-stable if there exists a series $F(n,x,v) = \sum_{k=0}^\infty \Phi_k(n,v)x^k$. Here $\Phi_k(n,v)$ are formal Laurent power series in $v$ with coefficients that are integer-valued quasi-polynomials in $n$.
We require that for all $k\in \mathbb{N}$ we have
\[
\lim_{n\to\infty}v^{-k(n+1)}\left(v^{-\delta(n)}f_n(v)-\sum_{j=0}^k\Phi_j(n,v)v^{j(n+1)} \right)=0
\]

In their paper they conjecture that $c$-stability holds for all Reshetikhin-Turaev knot invariants as one scales a dominant weight.

%British spelling/American spelling?
{\bf Acknowledgement.} We thank Stavros Garoufalidis for informing us about the notion of $c$-stability and the organisers of the 2018 Banff workshop on Modular forms and Quantum invariants for providing the right atmosphere for writing this paper. We thank Abhijit Champanerkar for interesting conversations about tails.
\bibliographystyle{plain}

\bibliography{References}

\end{document}